\documentclass{amsart}
\usepackage{amsmath}
\usepackage{amssymb}
\usepackage{amsthm}
\usepackage{epsfig}
\usepackage{amsfonts}
\usepackage{amscd}
\usepackage{mathrsfs}
\usepackage{enumerate}
\usepackage{latexsym}
\usepackage{url}

\newtheorem{theorem}{Theorem}[section]
\newtheorem{lemma}[theorem]{Lemma}

\newtheorem{notationandlemma}[theorem]{Notation and Lemma}
\newenvironment{proof*}{\noindent\textbf{Proof.}}{\hfill{\;$\blacksquare$}\newline}

\newenvironment{proof of the claim}{\noindent\textbf{Proof of the claim.}}{\hfill{$\square$}\newline}

\theoremstyle{definition}
\newtheorem{definition}[theorem]{Definition}

\newtheorem{example}[theorem]{Example}
\newtheorem{remark}[theorem]{Remark}

\numberwithin{equation}{section}

\begin{document}

\title[On closed subalgebras of $C_B(X)$]{On closed subalgebras of $C_B(X)$}

\author[M. Farhadi and M.R. Koushesh]{M. Farhadi and M.R. Koushesh$^*$}

\address{\textbf{[First author]} Department of Mathematical Sciences, Isfahan University of Technology, Isfahan 84156--83111, Iran.}

\email{m.farhadi@math.iut.ac.ir}

\address{\textbf{[Second author]} Department of Mathematical Sciences, Isfahan University of Technology, Isfahan 84156--83111, Iran and School of Mathematics, Institute for Research in Fundamental Sciences (IPM), P.O. Box: 19395--5746, Tehran, Iran.}

\email{koushesh@cc.iut.ac.ir}

\thanks{$^*$Corresponding author}

\thanks{This research was in part supported by a grant from IPM (No. 96030416).}



\begin{abstract}
For a completely regular space $X$ and a non-vanishing self-adjoint closed subalgebra $H$ of $C_B(X)$ which separates points from closed sets in $X$ we construct the Gelfand spectrum $\mathfrak{sp}(H)$ of $H$ as an open subspace of the compactification of $X$ generated by $H$. The simple construction of $\mathfrak{sp}(H)$ enables easier examination of its properties. We illustrate this by an example showing that the space $\mathfrak{sp}(H)$ is separable metrizable if and only if $H$ is countably generated.
\end{abstract}

\maketitle


\section{Introduction}

By a \textit{space} we mean a \textit{topological space}. We follow the definitions of \cite{E}, in particular, completely regular spaces and compact spaces are Hausdorff (consequently, locally compact spaces are completely regular), and Lindel\"{o}f spaces are regular. The field of scalars is the complex field $\mathbb{C}$, however, all results remain true (with exactly the same proof) in the real setting.

Let $X$ be a space. We denote by $C_B(X)$ the set of all continuous bounded scalar-valued mappings on $X$. The set $C_B(X)$ is a Banach algebra with pointwise addition and multiplication and the supremum norm. We denote by $C_0(X)$ the set of all $f$ in $C_B(X)$ which vanish at infinity (i.e., $|f|^{-1}([\epsilon,\infty))$ is compact for all $\epsilon>0$). For any $f$ in $C_B(X)$, the \textit{cozero-set} of $f$ is defined to be $X\setminus f^{-1}(0)$ and is denoted by $\mathrm{Coz}(f)$.

In \cite{FK}, assuming that $X$ is a completely regular space, we have represented a non-vanishing self-adjoint closed subalgebra $H$ of $C_B(X)$ which has local units as $C_0(Y)$ for some locally compact space $Y$. (Here, by $H$ having \textit{local units} we mean that for every closed subspace $C$ of $X$ and every neighborhood $U$ of $C$ in $X$ contained in the support of some element of $H$ there is some $h$ in $H$ which equals $1$ on $C$ and vanishes outside $U$.) Our purpose here is to improve our previous result by replacing the rather ``ad hoc'' assumption that $H$ has local units by the more standard assumption that $H$ separates points from closed sets in $X$. This will be accomplished by considering the compactification $\alpha_H X$ of $X$ generated by $H$ instead of the Stone--\v{C}ech compactification of $X$. Our approach here, apart from standardizing our assumption on $H$, enables us to derive ceratin properties of the space $Y$ (which coincide with the Gelfand spectrum of $H$). We illustrate this by an example in which we show that the space $Y$ is separable metrizable if and only if $H$ is countably generated. We conclude with the description of examples of completely regular spaces $X$ and non-vanishing self-adjoint closed subalgebras $H$ of $C_B(X)$ which separate points from closed sets in $X$, so that our conclusion may apply.

For other results of the same nature dealing with non-vanishing closed ideals of $C_B(X)$ (where $X$ is a completely regular space) we refer to \cite{K1}--\cite{K4}.

\section{The representation theorem}

We begin with the following definition.

\begin{definition}\label{PGS}
Let $X$ be a space. A subset $H$ of $C_B(X)$
\begin{itemize}
\item is \textit{self-adjoint} if $H$ contains the complex conjugate $\overline{h}$ of any element $h$ in $H$ (where $\overline{h}(x)=\overline{h(x)}$ for any $x$ in $X$).
\item is \textit{non-vanishing} if for any $x$ in $X$ there is some $h$ in $H$ such that $h(x)\neq 0$.
\item \textit{separates points in $X$} if for any distinct elements $x$ and $y$ in $X$ there is some $h$ in $H$ such that $h(x)\neq h(y)$.
\item \textit{separate points from closed sets in $X$} if for any closed subspace $C$ of $X$ and any $x$ in $X\setminus C$ there is some $h$ in $H$ such that $h(x)$ is not in $\mathrm{cl}_\mathbb{C}h(C)$.
\end{itemize}
\end{definition}

We need the following well known theorem. We include the proof here for notational convenience and completeness.

\begin{notationandlemma}\label{PJUG}
Let $X$ be a completely regular space and let $H$ be a subset of $C_B(X)$ which separates points from closed sets in $X$. For any $h$ in $H$, let $I_h=\mathrm{cl}_\mathbb{C}h(X)$. Let
\[e:X\longrightarrow\prod_{h\in H}I_h\]
be the evaluation map where
\[x\longmapsto\big\{h(x)\big\}_{h\in H}\]
for any $x$ in $X$. Then, $e$ is a homeomorphic embedding. Let
\[\alpha_H X=\mathrm{cl}_{\prod_{h\in H}I_h}e(X)\]
and identify $X$ with its image $e(X)$. Then $\alpha_H X$ is a compactification of $X$ on which every element of $H$ extends continuously. For any $h$ in $H$, we denote by $h_\alpha$ the (unique) continuous extension of $h$ on $\alpha_H X$.
\end{notationandlemma}

\begin{proof*}
It is known that $e$ is injective if $H$ separates points in $X$ and $e$ is a homeomorphic embedding if $H$ further separates points from closed sets in $X$. (See Theorem 2.3.20 of \cite{E}.) By our assumption $H$ separates points from closed sets in $X$ which also implies that $H$ separates points in $X$ (as all one point sets are closed in $X$). Note that $\prod_{h\in H}I_h$ is compact (as is a product of compact spaces) and therefore, its closed subspace $\alpha_H X$ is also compact. It is clear that $\alpha_H X$ contains $e(X)$ ($=X$) as a dense subspace. Therefore $\alpha_H X$ is a compactification of $X$. Note that $\pi_h(e(x))=h(x)$ for every $x$ in $X$, where $h$ is in $H$ and $\pi_h:\alpha_H X\rightarrow \mathbb{C}$ is the projection onto the $h$-th coordinate. Since $X$ is identified with $e(X)$, $\pi_h$ continuously extends $h$. Note that continuous extensions of $h$ coincide, as are identical on the dense subspace $e(X)$ ($=X$) of $\alpha_H X$.
\end{proof*}

In Theorem 2.9 of \cite{FK}, assuming that $X$ is a completely regular space, we have represented a non-vanishing self-adjoint closed subalgebra $H$ of $C_B(X)$ which has local units as $C_0(Y)$ for some locally compact space $Y$ (where as pointed previously, by $H$ having \textit{local units} it is meant that for every closed subspace $C$ of $X$ and every neighborhood $U$ of $C$ in $X$ contained in the support of some element of $H$ there is some $h$ in $H$ which equals $1$ on $C$ and vanishes outside $U$.) The space $Y$ has been constructed as an open subspace of the Stone--\v{C}ech compactification of $X$. The requirement on $H$ to have local units is rather ``ad hoc.'' In the next theorem, we improve Theorem 2.9 of \cite{FK} by replacing the assumption that $H$ has local units by the more standard assumption that $H$ separates points from closed sets in $X$. This will be done by considering the compactification $\alpha_H X$ of $X$ generated by $H$ instead of the Stone--\v{C}ech compactification of $X$.

In the proof of the next theorem we use the following corollary of the Stone--Weierstrass theorem. (See Chapter V of \cite{C}, Theorem 8.1 and Corollary 8.3.)

\begin{lemma}\label{PJDF}
Let $Y$ be a locally compact space and let $G$ be a non-vanishing self-adjoint closed subalgebra of $C_0(Y)$ which separates points in $Y$. Then $G=C_0(Y)$.
\end{lemma}

Also, we need to use the following version of the Banach--Stone theorem stating that the topology of a locally compact space $Y$ determines and is determined by the normed algebraic properties of $C_0(Y)$. (See Theorem 7.1 of \cite{B}. It turns out the even algebraic properties of $C_0(Y)$ suffice to determine the topology of $Y$; see \cite{AAN}.)

\begin{lemma}\label{KGSD}
For locally compact spaces $Y$ and $Z$, the normed algebras $C_0(Y)$ and $C_0(Z)$ are isometrically isomorphic if and only if the spaces $Y$ and $Z$ are homeomorphic.
\end{lemma}

We are now at a place to prove our representation theorem.

\begin{theorem}\label{FS}
Let $X$ be a completely regular space and let $H$ be a non-vanishing self-adjoint closed subalgebra of $C_B(X)$ which separates points from closed sets in $X$. Then $H$ is isometrically isomorphic to $C_0(Y)$ for some locally compact space
\[Y=\bigcup_{h\in H}\mathrm{Coz}(h_\alpha).\]
The space $Y$ is unique (up to homeomorphism) with this property. Moreover, the following are equivalent:
\begin{itemize}
\item[(1)] $H$ is unital.
\item[(2)] $H$ has $\mathbf{1}$.
\item[(3)] $Y$ is compact.
\item[(4)] $Y=\alpha_H X$.
\end{itemize}
\end{theorem}

\begin{proof*}
Let
\[\lambda_H X=\bigcup_{h\in H}\mathrm{Coz}(h_\alpha).\]
Note that $\lambda_H X$ contains $X$, as for any $x$ in $X$, since $H$ is non-vanishing, there is some $h$ in $H$ such that $h(x)\neq 0$, but $h_\alpha(x)=h(x)$. For any $h$ in $H$ we let
\[h_\lambda=h_\alpha|_{\lambda_H X}:\lambda_H X\longrightarrow \mathbb{C}.\]
Note that $h_\lambda$ extends $h$, as $h_\alpha$ does and $\lambda_H X$ contains $X$. We define the mapping
\[\phi:H\longrightarrow C_0(\lambda_H X)\]
by $\phi(h)=h_\lambda$ for any $h$ in $H$, and we check that $\phi$ is an isometric isomorphism.

Let $h$ be in $H$. For any $\epsilon>0$ we have
\[|h_\lambda|^{-1}\big([\epsilon,\infty)\big)=\lambda_H X\cap|h_\alpha|^{-1}\big([\epsilon,\infty)\big).\]
But the latter is $|h_\alpha|^{-1}([\epsilon,\infty))$, and is therefore, being closed in $\alpha_H X$, is compact. That is, $\phi(h)$ is in $C_0(\lambda_H X)$. That $\phi$ is a homomorphism and is injective follows from the observation that any two continuous mappings on $\lambda_H X$ which agree on the (dense) subspace $X$ of $\lambda_H X$ agree on the whole $\lambda_H X$. (As an example, $(f+g)_\lambda$ and $f_\lambda+g_\lambda$ coincide, as they are both identical to $f+g$ on $X$.) Note that for an $h$ in $H$, by continuity of $h$, we have
\[|h_\lambda|(\lambda_H X)=|h_\lambda|(\mathrm{cl}_{\lambda_H X}X)\subseteq\mathrm{cl}_\mathbb{C}|h_\lambda|(X)=\mathrm{cl}_\mathbb{C}|h|(X)\subseteq\big[0,\|h\|\big].\]
Thus $\|h_\lambda\|\leq\|h\|$. That $\|h\|\leq\|h_\lambda\|$ is clear, as $h_\lambda$ extends $h$. Therefore $\|h_\lambda\|=\|h\|$. This shows that $\phi$ preservers norm. To conclude the proof we need to show that $\phi$ is surjective. By Lemma \ref{PJDF} we suffice to check that $\phi(H)$ is a non-vanishing self-adjoint closed subalgebra of $C_0(\lambda_H X)$ which separates points in $\lambda_H X$.

As it is observed in the above lines $\phi(H)$ is a subalgebra of $C_0(\lambda_H X)$, which is also closed in $C_0(\lambda_H X)$, as is an isometrically isomorphic image of a complete normed space. The fact that $\phi(H)$ is non-vanishing on $\lambda_H X$ follows from the definition of $\lambda_H X$. (Indeed, for any $z$ in $\lambda_H X$ we have $h_\alpha(z)\neq 0$, and in particular $h_\lambda(z)\neq 0$, for some $h$ in $H$.) It is also clear that $\phi(H)$ is self-adjoint, as $\overline{(h_\lambda)}=(\overline{h})_\lambda$ for any $h$ in $H$ (since the two mapping coincide with $\overline{h}$ on the dense subspace $X$ of $\lambda_H X$). We check that $\phi(H)$ separates points in $\lambda_H X$. We suffice to check that $\{h_\alpha:h\in H\}$ separates points in $\alpha X$. But this follows, as if $s=\{s_h\}_{h\in H}$ and $t=\{t_h\}_{h\in H}$ are distinct elements of $\alpha X$, then $s_h\neq t_h$ for some $h$ in $H$, and then by Lemma \ref{PJUG} (and its proof), we have
\[h_\alpha(s)=\pi_h(s)=s_h\neq t_h=\pi_h(t)=h_\alpha(t).\]
Note that $\lambda_H X$ is locally compact, as (by its definition) is open in the compact space $\alpha_H X$. Lemma \ref{PJDF} now implies that
\[\phi(H)=C_0(\lambda_H X).\]
This shows that $H$ is isometrically isomorphic to $C_0(\lambda_H X)$. The uniqueness assertion of the theorem follows from Lemma \ref{KGSD}.

Finally, we check that (1)--(4) are equivalent. To see that (1) implies (2), let $u$ be a unit in $H$. For any $x$ in $X$ let $h_x$ be an element of $H$ which does not vanish at $x$. Then $u(x)h_x(x)=h_x(x)$ which implies that $u(x)=1$. That (2) implies (4) is clear and follows from the definition of $\lambda_H X$. That (2) implies (1), and (4) implies (3) are also clear. Finally, (3) implies (1), as if $Y$ is compact then $C_0(Y)$ coincides with $C_B(Y)$ and is therefore unital (as $C_B(Y)$ is). But then $H$ is unital, as is isometrically isomorphic to $C_0(Y)$ by the first part.
\end{proof*}

\begin{remark}\label{JBB}
In Theorem \ref{FS}, for a completely regular space $X$ and a non-vanishing self-adjoint closed subalgebra $H$ of $C_B(X)$ which separates points from closed sets in $X$ we have proved that $H$ and $C_0(\lambda_HX)$ are isometrically isomorphic. On the other hand (in the case when the field of scalars is $\mathbb{C}$) by the commutative Gelfand--Naimark theorem, $H$, being a Banach algebra, is isometrically isomorphic to $C_0(Y)$ with $Y$ being the \textit{spectrum} (or the \textit{maximal ideal space}) of $H$ with the Gelfand (or Zariski) topology. The uniqueness part of Theorem \ref{FS} implies that $\lambda_HX$ and the spectrum of $H$ coincide. Indeed, our result here has the advantage that it explicitly constructs the spectrum of $H$ as a subspace of a compactification of $X$ with a known structure. As we will see next, this may provide some information which is not generally expected to be derivable from the standard Gelfand theory. For this purpose, Theorem \ref{OUJH} may be considered as a particular example, in which we prove that the spectrum of $H$ is separable metrizable if $H$ is countably generated.
\end{remark}

The above remark motivates to introduce the following notation.

\begin{definition}\label{PGS}
Let $X$ be a completely regular space and let $H$ be a non-vanishing self-adjoint closed subalgebra of $C_B(X)$ which separates points from closed sets in $X$. We denote by $\mathfrak{sp}(H)$ the unique locally compact space $Y$ such that $H$ and $C_0(Y)$ are isometrically isomorphic.
\end{definition}

Our next theorem should be known to various authors; we deduce it here, however, as an easy corollary of our representation theorem. We need to use the following lemma.

\begin{lemma}\label{OJH}
Let $X$ be a completely regular space and let $H$ be a subset of $C_B(X)$ which separates points from closed sets in $X$. Then
\[\alpha_{\overline{H}}X=\alpha_HX.\]
Here the bar denotes the closure in $C_B(X)$.
\end{lemma}

\begin{proof*}
Note that $\overline{H}$ separates points from closed sets in $X$, as $H$ does. By the construction in Lemma \ref{PJUG} we have
\[X\subseteq\alpha_HX\subseteq\alpha_{\overline{H}}X,\]
with $X$ being dense in the latter. This implies that $\alpha_{\overline{H}}X\subseteq\alpha_HX$ which together with the above relation proves the lemma.
\end{proof*}

In the following theorem we use the known fact that a locally compact metrizable space $Y$ is $\sigma$-compact if and only if $C_0(Y)$ is separable. (See Theorem 3.5.17 of \cite{C0}.) Observe that in locally compact metrizable spaces $\sigma$-compactness and separability coincide. (See Corollary 4.1.16 and Exercise 3.8.C(b) of \cite{E}.)

For a space $X$ and a subset $G$ of $C_B(X)$ we define the algebra (closed algebra, respectively) generated by $G$ as the smallest subalgebra (closed subalgebra, respectively) of $C_B(X)$ which contains $G$. We denote by $\langle G\rangle$ the subalgebra of $C_B(X)$ generated by $G$. Note that the closed subalgebra of $C_B(X)$ generated by $G$ is then the closure $\overline{\langle G\rangle}$ of $G$ in $C_B(X)$.

\begin{theorem}\label{OUJH}
Let $X$ be a completely regular space and let $H$ be a non-vanishing self-adjoint closed subalgebra of $C_B(X)$ which separates points from closed sets in $X$. Then $H$ is countably generated if and only if $\mathfrak{sp}(H)$ is separable and metrizable.
\end{theorem}

\begin{proof*}
Suppose that $\mathfrak{sp}(H)$ is a separable metrizable space. Since $\mathfrak{sp}(H)$ is locally compact, $\mathfrak{sp}(H)$ is then $\sigma$-compact, and therefore $C_0(\mathfrak{sp}(H))$ is separable. But then $H$ is also separable (as is isometrically isomorphic to $C_0(\mathfrak{sp}(H))$) and therefore countably generated.

Suppose that $H$ is countably generated. Let
\[H=\overline{\langle h_1,h_2,\ldots\rangle},\]
where $h_1,h_2,\ldots$ are in $H$ and the bar denotes the closure in $C_B(X)$. Note that
\[\langle h_1,h_2,\ldots\rangle=\Big\{\sum_{i\in I} c_i h^{p_1}_{k_1}\cdots h^{p_{n_i}}_{k_{n_i}}:\mbox{$I$ is finite, $c_i\in\mathbb{C}$, $k_j,p_j\in\mathbb{N}$ for $j=1,\ldots,n_i$}\Big\}.\]
Let
\[Q=\Big\{\sum_{i\in I} q_i h^{p_1}_{k_1}\cdots h^{p_{n_i}}_{k_{n_i}}:\mbox{$I$ is finite, $q_i\in\mathbb{Q}\times\mathbb{Q}$, $k_j,p_j\in\mathbb{N}$ for $j=1,\ldots,n_i$}\Big\}.\]
One can check that $Q$ is dense in $\langle h_1,h_2,\ldots\rangle$, and is therefore dense in $H$, i.e., $H=\overline{Q}$. Also, $Q$ separates points from closed sets in $X$. (To see this, let $C$ be a closed subspace of $X$ and let $x$ be an element of $X\setminus C$. Then there is an element $h$ in $H$ such that $h(x)$ is not in $\mathrm{cl}_\mathbb{C}h(C)$. Let $\epsilon>0$ such that $|h(x)-h(c)|>\epsilon$ for every $c$ in $C$. Let $f$ be an element of $Q$ such that $\|f-h\|<\epsilon/3$. Then
\[\big|f(x)-f(c)\big|\geq\big|h(x)-h(c)\big|-\big|f(c)-h(c)\big|-\big|f(x)-h(x)\big|>\epsilon-\epsilon/3-\epsilon/3=\epsilon/3\]
for every $c$ in $C$, and therefore $f(x)$ is not in $\mathrm{cl}_\mathbb{C}f(C)$.) By Lemma \ref{OJH} we have
\[\alpha_HX=\alpha_QX.\]
Note that $Q$ is countable. Therefore $\alpha_QX$ is a separable metrizable space, as $\alpha_QX$ (by the construction in Lemma \ref{PJUG}) will be a subspace of a countable product of spaces $\mathbb{C}$. But then $\lambda_HX$ is also a separable metrizable space, as $\lambda_HX$ is a subspace of $\alpha_HX$ ($=\alpha_QX$). Finally, observe that $\lambda_HX$ is the spectrum of $H$ by Theorem \ref{FS} (and its proof).
\end{proof*}

In the following we give examples of completely regular spaces $X$ and non-vanishing self-adjoint closed subalgebras $H$ of $C_B(X)$ which separate points from closed sets in $X$, thus, satisfying the assumption in Theorem \ref{FS}.

\begin{example}\label{JH}
Let $\mathscr{P}$ be a topological property which is closed hereditary (in the sense that every closed subspace of a space with $\mathscr{P}$ has $\mathscr{P}$) and preserved under countable closed unions of subspaces (in the sense that every space which is a countable union of closed subspaces with $\mathscr{P}$ has $\mathscr{P}$). (Examples of such topological properties are the Lindel\"{o}f property, subparacompactness, submetacompactness, the submeta-Lindel\"{o}f property, weakly $\theta$-refinability and weakly $\delta\theta$-refinability; see Theorems 7.1 and 7.3 of \cite{Bu}.) Let $X$ be a completely regular locally-$\mathscr{P}$ space (in the sense that every element of $X$ has a neighborhood in $X$ with $\mathscr{P}$). Let
\[H=\big\{f\in C_B(X):\mbox{$|f|^{-1}\big([\epsilon,\infty)\big)$ has $\mathscr{P}$ for all $\epsilon>0$}\big\}.\]
Then $H$ is a non-vanishing self-adjoint closed subalgebras of $C_B(X)$ which separates points from closed sets in $X$, as we now check. (In Theorem 3.3.9 of \cite{K4} it has been shown that $H$ is a non-vanishing closed ideal of $C_B(X)$.) That $H$ is closed under addition and multiplication follows from the fact that for any two elements $f$ and $g$ in $H$ and $\epsilon>0$ we have
\[|f+g|^{-1}\big([\epsilon,\infty)\big)\subseteq|f|^{-1}\big([\epsilon/2,\infty)\big)\cup|g|^{-1}\big([\epsilon/2,\infty)\big),\]
where the latter (being the union of two closed subspaces with $\mathscr{P}$) has $\mathscr{P}$, and
\[|fg|^{-1}\big([\epsilon,\infty)\big)\subseteq|f|^{-1}\big([\epsilon/M,\infty)\big)\]
where $M>\|g\|$. Therefore, $|fg|^{-1}([\epsilon,\infty))$ and $|f+g|^{-1}([\epsilon,\infty))$, being closed subspaces of spaces with $\mathscr{P}$, have $\mathscr{P}$. Similarly, one can check that $H$ is closed under scalar multiplication. Therefore, $H$ is a subalgebra of $C_B(X)$. To check that $H$ is a closed subalgebra of $C_B(X)$, let $h_n\rightarrow f$, where $h_1,h_2,\ldots$ is a sequence in $H$ and $f$ in $C_B(X)$. Then, for every $\epsilon>0$ we have
\[|f|^{-1}\big([\epsilon,\infty)\big)\subseteq\bigcup_{n=1}^\infty|h_n|^{-1}\big([\epsilon/2,\infty)\big),\]
where the latter (being a countable union of closed subspaces with $\mathscr{P}$) has $\mathscr{P}$, and therefore, so does its closed subspace $|f|^{-1}([\epsilon,\infty))$. Therefore, $f$ is in $H$, and thus $H$ is closed in $C_B(X)$. It is clear that $H$ is self-adjoint. To check that $H$ is non-vanishing and separates points from closed sets in $X$, let $C$ be a closed subspace of $X$ and let $x$ be in $X\setminus C$. Since $X$ is locally-$\mathscr{P}$, there is a neighborhood $U$ of $x$ in $X$ which has $\mathscr{P}$. Let $W=U\cap(X\setminus C)$. Then $W$ is a neighborhood of $x$ in $X$. By complete regularity of $X$ there is some $f:X\rightarrow[0,1]$ such that $f(x)=1$ and $f(y)=0$ for all $y$ outside $W$, in particular, $f(C)=0$. Note that $f$ is in $H$, as $|f|^{-1}([\epsilon,\infty))\subseteq U$ for every $\epsilon>0$, $U$ has $\mathscr{P}$, and $|f|^{-1}([\epsilon,\infty))$ is closed in $U$.
\end{example}

\end{document}